\documentclass[11pt]{article}
\usepackage{amsmath}
\usepackage{amsthm,amsfonts,amssymb,latexsym}
\usepackage[utf8]{inputenc}
\usepackage[english]{babel}
\usepackage{color}

\usepackage{hyperref} 
\hypersetup{
colorlinks=true, 
breaklinks=true, 
urlcolor= blue, 
linkcolor= blue, 
}
\newtheorem{theorem}{Theorem}[section]
\newtheorem{prop}[theorem]{Proposition}
\newtheorem{lem}[theorem]{Lemma}
\newtheorem{fact}[theorem]{Fact}
\newtheorem{cor}{Corollary}[theorem]

\theoremstyle{definition}
\newtheorem{defn}{Definition}
\theoremstyle{remark}
\newtheorem{exe}{Example}

\def\R{{\mathbb R}}
\def\N{{\mathbb N}} 
\newcommand{\Rex}{\R \cup \{\infty\}} 
\newcommand{\imp}{\Rightarrow}
\newcommand{\pmi}{\Leftarrow}
\newcommand{\tq}{:}
\newcommand{\la}{\langle}
\newcommand{\ra}{\rangle}

\newcommand{\lsc}{lsc}
\newcommand {\dom} {{\rm dom} \kern.15em} 
\newcommand {\diam} {{\rm diam} \kern.15em} 
\newcommand {\epi} { {\rm epi}\kern .15em} 
\newcommand {\hypo} { {\rm hypo}\kern .15em} 
\newcommand {\graph} {{\rm graph} \kern .15em}
\newcommand {\Li} {\mathrm{Li} \kern .15em}
\newcommand {\Ls} {\mathrm{Ls} \kern .15em}
 
\newcommand{\gap}{{\rm D}}

\newcommand{\lip}{{\rm lip}\kern.15em}
\newcommand{\cl}{{\rm cl}\kern.15em}
\newcommand {\inT} {{\rm int} \kern.15em} 
\newcommand{\del} {\partial}
\newcommand {\dis} {\displaystyle} 
\newcommand {\eps} {\varepsilon} 
\newcommand {\ldx} {\lambda_x} 
\newcommand {\ld} {\lambda} 
\newcommand {\bx} {\bar{x}} 
\newcommand {\xb} {\bar{x}} 
\newcommand {\yb} {\bar{y}} 
\newcommand {\zb} {\bar{z}} 
\newcommand {\xbs} {\bar{x}^*} 
\newcommand {\ybs} {\bar{y}^*} 
 
\newcommand {\ind}{\delta}
\newcommand {\LSC} {{\rm LSC}} 
\newcommand{\pf}{\noindent{\bf Proof.~}} 
\newcommand{\cqfd}{\mbox{}\nolinebreak\hfill\rule{2mm}{2mm}\bigbreak}
\newcommand{\finpf}{\cqfd}
\newcommand{\lda}{\lambda}
\newcommand{\deltb}{\bar{\delta}}
\newcommand{\ldab}{\bar{\lambda}}
\newcommand{\PK}{{Kuratowski }}

\newcommand{\V}{\mathcal{N}}

\begin{document} 
\begin{center}
{\large\bf\sc
Subdifferential stability and subdifferential sum rules
}
\medskip\\
\today
\end{center}

\begin{center}
  {\small\begin{tabular}{c}
  Marc Lassonde\\
  LIMOS, Clermont-Ferrand, France\\
  E-mail: marc.lassonde@gmail.com
  \end{tabular}}
\end{center}
\begin{center}
{\textit{Dedicated to Alex Ioffe on the occasion of his 80th birthday}}
\end{center}

\medbreak\noindent
\textbf{Abstract.}
In the first part, we discuss the stability of the
strong slope and of the subdifferential of a
lower semicontinuous function with respect
to Wijsman perturbations of the function, i.e.\ perturbations
described via Wijsman convergence.
In the second part, we show how subdifferential sum rules can be
viewed as special cases of subdifferential stability results.

\medbreak\noindent
\textbf{Keywords:}
  lower semicontinuity, Wijsman convergence, slope,
  subdifferential calculus, trustworthiness.
  
\medbreak\noindent
\textbf{2010 Mathematics  Subject Classification:}
  49J52, 49J53, 49J45, 26E15.

\section{Introduction}
In the first part, we study the stability of the
strong slope \cite{DMT80} and of the subdifferential of a
lower semicontinuous function with respect to variational perturbations
of the function. This issue was initiated in our work \cite{GL03}.
There we showed that the slope of the sum $f+x^*$ of
a lower semicontinuous function $f$ and of a continuous linear functional $x^*$
is stable under \textit{slice} perturbations of $f$,
where the notion of slice convergence for lower semicontinuous,
non necessarily convex, functions was introduced
in \cite{GL00}.
Here we show instead that the slope of $f$ is stable under \textit{Wijsman}
perturbations of $f$, a weaker and classical notion of convergence
(see \cite{Bee93,Wij66}).
Our preceding result can be recovered from the following fact
(see Theorem \ref{W-versus-slice}):
{\em A sequence $(f_n)_n$ is slice
convergent to $f$ if and only if the sequence
$(f_n+x^*)_n$ is Wijsman convergent to $f+x^*$ for every $x^*\in X^*$}.
Applications to the stability of \textit{trustworthy} subdifferentials
(see \cite{Iof12}) under Wijsman convergence are adapted from \cite{GL03}. 

In the second part, we show how subdifferential sum rules can be
viewed as special cases of subdifferential stability results.

The results of this article were largely announced in \cite{Las02}.
\medbreak
{\it Notation.}
Except where otherwise stated, $X$ stands for a real Banach space
and $X^\ast$ for its topological dual. All functions are assumed
to be extended-real-valued and lower semicontinuous (\lsc);
$\LSC(X)$ denotes the space of all such functions on $X$.
For $f\in\LSC(X)$, we denote by $\dom f: =\{x\in X\tq f(x)<\infty\}$
the \textit{effective domain} of $f$, by 
$\graph f :=\{\, (x,\alpha)\tq f(x)=\alpha\,\}$ the \textit{graph }of $f$,
by 
$\epi f :=\{\, (x,\alpha)\tq f(x)\le\alpha\,\}$ the \textit{epigraph} of $f$
and by
$\hypo f :=\{\, (x,\alpha)\tq f(x)\ge\alpha\,\}$ the \textit{hypograph} of $f$.
We write $x\to_f\xb$ to say that $x\to\xb$ and $f(x)\to f(\xb)$.
For any two functions $f,g\in\LSC(X)$ we denote by
$$f\triangledown g:x\mapsto (f\triangledown g)(x):=\inf_{z\in X} (f(z)+g(x-z))$$
the \textit{inf-convolution} of $f$ and $g$.
The closed $\eps$-ball centered at point $x$ is written 
$B_\eps(x)$.
For a subset $S\subset X$ and a norm $\|\cdot\|$ on $X$,
the distance of a point $x\in X$ to $S$ is given by
$$
d_S(x)=d(x,S):=\inf\,\{\,\|x-a\|\tq a\in S\,\}\,,
$$ 
and the closed uniform $\delta$-neighborhood of $S$ ($\delta\ge 0$) is
defined by
$$
B_\delta(S):=\{x\in X \tq d_S(x)\le\delta\}.
$$
The diameter of $S$ is given by
$
\diam(S) := \sup\,\{\,\|x-y\|\tq x,y\in S\,\}\,,
$
and the indicator of $S$ is the function $\delta_S:X\to\Rex$ defined by
$$
\delta_S(x) :=
\left \{
	  \begin{array}{ll}
	  0 & \mbox{if }~x\in S\\
	  \infty & \mbox{otherwise}.
	  \end{array}
\right .
$$
For $f:X\to\Rex$ et $S\subset X$, we write
$f_S:=f+\delta_S$ for the `restriction' of $f$ to $S$,
and $\inf_S f:=\inf f(S)$.

\section{Convergences of sets}
\label{setconv}
We recall (see \cite[Definition 5.2.1]{Bee93}) that
the {\em lower\/} and {\em upper limits\/} of a sequence of
sets $(S_n)_n$ in a Hausdorff topological space $Y$ are respectively
defined by
\begin{align}
\Li S_n&:=\{\, y\in Y\tq \forall V\,\in\V(y),\exists N\in \N,
\, \forall n\geq N:\, S_n\cap V\not=\emptyset\,\},\label{liminf}
\medskip\\
\Ls S_n&:=\{\, y\in Y\tq \forall \,V\in \V(y),\, \forall  N\in \N,
\, \exists n\ge N :\, S_n\cap V\not=\emptyset\,\}.\label{limsup}
\end{align}
In a metric space $Y$ these formulas reduce to
\begin{align}
\Li S_n&=\{\, y\in Y\tq \lim_n d(y,S_n)=0\},\tag{\ref{liminf}b}\label{liminfb}
\medskip\\
\Ls S_n&=
\{\, y\in Y\tq \liminf_n d(y,S_n)=0\}.\tag{\ref{limsup}b}\label{limsupb}
\end{align}
Definitions \eqref{liminfb}--\eqref{limsupb} go back to Peano (1887, 1908),
definitions \eqref{liminf}--\eqref{limsup} were popularized by Kuratowski (1948):
see Dolecki-Greco \cite{DG07} for historical comments.

Lower and upper limits of a sequence of sets
describe two symmetric
behaviors of the sequence with respect to individual points of the space:
a point $y$ is in $\Li  S_n$ if and only if every neighborhood
of $y$ `hits' $S_n$ eventually, while
a point $y$ is not in $\Ls  S_n$ if and only if some neighborhood
of $y$ `misses' $S_n$ eventually.
In normed spaces, this hit-and-miss behavior is more conveniently
described by using gap distances. We recall that
{\em the gap distance\/}
between two sets $A$, $B$ of $(Y,\|.\|)$ is given by 
$$ 
\gap(A,B) := \inf \{\, \|a-b\| \tq a\in A, \; b\in B\,\}
          = \inf \{\, d(a,B) \tq a\in A\,\}.
$$ 

\begin{prop}\label{hit-miss-lim}
Let $Y$ be a normed space and let $\{S_n, S\}\subset Y$
with $S$ closed. Then: 
\smallbreak
{\rm (a)} $\begin{array}[t]{lcl}
\dis S\subset \Li S_n &\Longleftrightarrow&
\dis \forall y\in Y, \quad\limsup_n\, d(y,S_n)\le d(y,S)\\
&\Longleftrightarrow&
\dis \forall S'\subset Y, \quad\limsup_n\, \gap(S',S_n)\le \gap(S',S).
\end{array}$
\smallbreak
{\rm (b)} $\dis ~S\supset \Ls S_n ~\Longleftrightarrow\smallskip\\
\hspace*{\fill}\forall y\in Y, \quad d(y,S)>0\imp
\sup_{\delta>0}\liminf_n \gap(B_\delta(y),S_n)>0.$
\smallbreak
{\rm (c)}
$\dis ~d(y,S) \leq \liminf_n d(y,S_n)~\Longleftrightarrow\smallskip\\
\hspace*{\fill}\forall \lda\ge 0, \quad
\gap(B_\lda(y),S)>0 \imp
           \sup_{\delta>0}\liminf_n \gap(B_{\lda+\delta}(y),S_n) >0.$
\end{prop}
\pf
The above assertions are certainly well known, but we provide a full
proof for completeness.

(a) The first property implies the second one: if
$\gamma>d(y,S)$,
then the open set $V:=\inT B(y,\gamma)$ contains a point $z$
in $S\subset \Li S_n$, and
since $V\in\V(z)$, it follows from \eqref{liminf} that
$S_n\cap V\neq\emptyset$ eventually,
hence $\limsup_n d(y,S_n) \le \gamma$, as required.
The second property implies the third one: if
$\gamma>\gap(S',S)$, then $d(z,S)<\gamma$ for some $z\in S'$, so
$\limsup_n \gap(S',S_n)\leq \limsup_n d(z,S_n)\leq d(z,S)<\gamma$.
The third property implies the first one: if $y\in S$, then
$\limsup_n d(y,S_n)=0$ by the third property, so $y\in \Li S_n$
by \eqref{liminfb}.

\smallbreak
(b) To prove `$\imp$', let $d(y,S)>0$.
Then, $y$ is not in $S\supset\Ls S_n$, so,
according to the definition \eqref{limsupb}
of the upper limit, there exist $\delta>0$ and
$N\in\N$ such that $d(y,S_n)>\delta$ for all $n\ge N$.
Hence, for any $\delta'\in (0,\delta)$, we have
$\gap(B_{\delta'}(y),S_n)\ge \delta-\delta'>0$, showing that the
property in the second half of (b) holds.
To prove `$\pmi$', let $y\notin S$. Then,
$d(y,S)>0$, so, for some $\delta>0$,
$\gap( B_\delta(y),S_n) >0$ eventually,
which implies that $S_n\cap  B_\delta(y)=\emptyset$ eventually, that is,
$y\notin\Ls S_n$.

\smallbreak
(c) To prove `$\imp$', let $\gamma:=\gap(B_\lda(y),S)>0$. Then,
$d(y,S)\ge \lda+\gamma$:
indeed, if $z\in S$, the point $y'=y+\lda(z-y)/\|z-y\|$ is in
$B_\lda(y)$, so $\|z-y\|=\lda+\|z-y'\|\ge\lda+\gamma$. It follows
from our assumption that
$d(y,S_n)\ge \lda+\gamma$ eventually, hence,
as easily seen,
for any $\delta\in(0,\gamma)$ we have
$\gap(B_{\lda+\delta}(y),S_n)\ge\gamma-\delta>0$ eventually,
showing that the second property in (c) holds.
To prove `$\pmi$', assume $d(y,S)>\lda$.
Then, $\gap(B_\lda(y),S)>0$, so, $\gap( B_{\lda}(y),S_n) >0$ eventually,
hence $d(y,S_n)>\lda$ eventually, as required.
\finpf

A sequence of sets $(S_n)_n$ is declared {\em Kuratowski convergent\/} to $S$
(or {\em Peano-Kuratowski convergent\/} to $S$)
provided $\Li S_n=\Ls S_n= S$ \cite[Definition 5.2.3]{Bee93}, whereas
in a metric space $(Y,d)$,
the sequence is declared {\em Wijsman convergent\/}
to $S$ provided $d(y, S_n)\to d(y,S)$ for every $y\in Y$
\cite[Definition 5.2.3]{Bee93}.
It readily follows from Proposition~\ref{hit-miss-lim}
that, in normed spaces, both convergences can be
characterized by a hit-and-miss criterion using gap distances:

\begin{cor}\label{hit-miss}
Let $(Y,\|.\|)$ be a normed space and let $\{S_n, S\}\subset Y$
with $S$ closed. Then:
\smallbreak
{\rm (a)}
$(S_n)_n$ is Wijsman convergent to $S$
if and only if, for every $y\in Y$ and $\lda\ge 0$,
\medskip\\
{\rm (W)}$ \left\lbrace
    \begin{array}{rl} 
  {\rm (i)}  &y\in S \Rightarrow\lim_n d(y,S_n)=0;\smallbreak\\ 
  {\rm (ii)}  & \gap(B_\lda(y),S)>0 \Rightarrow 
           \sup_{\delta>0}\liminf_n \gap(B_{\lda+\delta}(y),S_n) >0.
    \end{array} 
\right. $
\medbreak
{\rm (b)}
$(S_n)_n$ is Kuratowski convergent to $S$ if and only if ${\rm (W)}$
holds for every $y\in Y$ and $\lda= 0$.
\end{cor}

Since  $\Li S_n$ and $\Ls S_n$ are always closed sets, the limit $S$ of
a Kuratowski or Wijsman convergent sequence is always a closed set. 

\section{Convergences of functions}
\label{funcconv}
A sequence of functions $(f_n)_n\subset \LSC(X)$
is declared {\em epi-convergent\/} (or {\em $\Gamma$-convergent\/}) to 
a function $f\in\LSC(X)$
provided the sequence of their epigraphs $(\epi f_n)_n$ is \PK convergent
to $\epi f$ in $X\times \R$ \cite[Definition 5.3.1]{Bee93}.
Likewise, the sequence $(f_n)_n$
is declared {\em Wijsman convergent\/} to $f$ provided the sequence
$(\epi f_n)_n$ is Wijsman convergent to $\epi f$ in $X\times \R$
supplied with the max norm $\|(x,t)\| := \max\{\|x\|,|t|\}$.

In view of Corollary~\ref{hit-miss}, epi-convergence and Wijsman convergence
of functions are characterized by hit-and-miss criteria. The proof uses
the following easy observation:
\begin{fact}[{\cite[Lemma 2.2]{GL00}}]\label{epi-hypo-lemma}
For any two functions $f,g:X\to\Rex$ one has\smallskip\\
\centerline{$
\gap(\hypo g,\epi f)=\gap(\hypo g,\graph f)=
\gap(\graph g,\epi f).$}
\end{fact}

\begin{prop}\label{hit-miss-varcv}
Let $X$ be a normed space and $\{f_n,f\}\subset \LSC(X)$.
Then:\smallbreak
{\rm (a)} $(f_n)_n$ is Wijsman convergent to $f$ if and only if
for every $x\in X$, $\lda\ge 0$ and $\alpha\in\R$,
\smallskip\\
$  {\rm (W_g)} \left\lbrace
    \begin{array}{rl} 
  {\rm (i)}  &\exists\, (x_n)_n\subset X: \;  
  x_n \to x, \; f_n(x_n) \to f(x);\smallbreak\\ 
  {\rm (ii)}  & \dis\gap(B_\lda(x)\times\alpha,\epi f)>0 \Rightarrow
\sup_{\delta>0}\liminf_n \gap(B_{\lda+\delta}(x)\times\alpha,\epi f_n) >0.
    \end{array} 
\right. $ 
\medbreak
{\rm (b)} $(f_n)_n$ is epi-convergent to $f$ if and only if ${\rm (W_g)}$
holds for every $x\in X$, $\lda= 0$ and $\alpha\in\R$. 
\end{prop}
\pf
We have to show that the formulas ${\rm (W_g)}$ are equivalent to
the corresponding formulas (W) in Corollary~\ref{hit-miss}
with $S_n$, $S$ replaced
by $\epi f_n$, $\epi f$ and
$B_\lda(y)$ replaced by the balls of $X\times \R$ supplied with
the box norm, namely $B_\lda(x)\times [\beta,\alpha]$. 

We first observe that (W)(ii) and ${\rm (W_g)}$(ii) are equivalent:
this is due to the
fact that for any $h$, $\gap(B_\lda(x)\times[\beta,\alpha],\epi h)=
\gap(B_\lda(x)\times\alpha,\epi h)$, as
follows from Fact \ref{epi-hypo-lemma}
with $g:=\alpha +\delta_{B_\lda(x)}$:
\begin{eqnarray*}
\gap(\hypo g,\epi h) &\le& \gap(B_\lda(x)\times[\beta,\alpha],\epi h)\\
                &\le& \gap(B_\lda(x)\times\alpha,\epi h)\\
                & =& \gap(\graph g,\epi h)\\
                & =& \gap(\hypo g,\epi h).
\end{eqnarray*}

Next, we claim that (W)(i) is equivalent to
$${\rm (i')}~~ x\in\dom f\imp \exists\, (x_n)_n\subset X: \;  
  x_n \to x,\,\limsup_n f_n(x_n) \le f(x).$$
Indeed, assume (W)(i) and let 
$x\in\dom f$. It follows from (W)(i) that
$\lim_n d((x,f(x)),\epi f_n)=0$, hence 
there exists a sequence $((x_n,\alpha_n))_n$ in $\epi f_n$
such that $x_n \to x$ and $\alpha_n \to f(x)$. This implies that
$\limsup_n f_n(x_n) \leq \limsup_n \alpha_n = f(x)$.
Conversely, assume that (i') holds, let $(x,\alpha)\in\epi f$ and let
$\eps>0$. It follows from (i') that for some sequence $x_n\to x$,
$f_n(x_n)<\alpha+\eps$ eventually. We therefore have 
$$\limsup_n d((x,\alpha),\epi f_n)\le
\limsup_n d((x,\alpha),(x_n,\alpha+\eps))=\eps.$$
This shows that $\lim_n d((x,\alpha),\epi f_n)=0$, as required.

To complete the proof of the proposition,
it suffices to show that ${\rm (W_g)}$(ii) implies that
$f(x)\le\liminf_n f_n(x_n)$ for every sequence $x_n\to x$. Let
$\alpha< f(x)$. Then, $d((x,\alpha),\epi f)>0$, so by
${\rm (W_g)}$(ii) with $\lda=0$, for some $\delta>0$ one has
$\gap(B_{\delta}(x)\times\alpha,\epi f_n) >0$ eventually. Now, let
$x_n\to x$. Then $x_n\in B_{\delta}(x)$ eventually, so 
$d((x_n,\alpha),\epi f_n) >0$ eventually, i.e., $\alpha<f_n(x_n)$
eventually, which was to be proved.
\finpf

Formulas ${\rm (W_g)}$ suggest possible localizations of
either concept of variational convergence `at a given point $x$'; we
consider only Wijsman convergence:

\begin{defn}
A sequence of functions $(f_n)_n\subset\LSC(X)$ is
{\em Wijsman convergent to $f$ at $x$ with radius $\lda_x\in(0,+\infty]$\/}
provided ${\rm (W_g)}$ holds at $x$ for all
$\lda\in [0,\lda_x)$ and all $\alpha\in\R$.
\end{defn}
Evidently, (global) Wijsman convergence implies
(local) Wijsman convergence at every point. The converse
need not be true.

\smallbreak
In \cite{GL00,GL03}, a stronger concept of variational convergence
is considered: roughly, it consists in demanding ${\rm (W_g)}$ to hold
not only for horizontal bounded slices $B_\lda(x)\times\alpha$, i.e.\
graphs of constant maps restricted to balls, but
more generally for all non-vertical bounded slices, i.e.\ graphs of
continuous affine maps restricted to balls. The localization of this
concept of convergence at an individual point reads as follows:
\begin{defn}\label{variationalconv-slice}
A sequence of functions $(f_n)_n\subset \LSC(X)$ is {\em slice convergent\/}
{\em to $f$ at $x\in X$ with radius $\lda_x\in(0,+\infty]$\/}
if for all functions $\varphi_\lda:=\varphi+\ind_{B_\lda(x)}$, with
$\varphi$ affine continuous and $\lda\in[0,\ldx)$,
\medskip\\
$ {\rm (s)} \left\lbrace \!\!
    \begin{array}{rl} 
  {\rm (i)}  & \exists\, (x_n)_n\subset X: \;  
  x_n \to x, \; f_n(x_n) \to f(x);\smallbreak\\ 
  {\rm (ii)}  & \dis\gap(\graph \varphi_\lda,\epi f)>0 \Rightarrow
\sup_{\delta>0}\liminf_n  \gap(\graph \varphi_{\lda+\delta},\epi f_n)>0. 
    \end{array} 
\right. $ 
\medskip\\
The sequence $(f_n)_n$ is (globally) {\em slice convergent\/} to $f$ if
{\rm (s)} holds at every point $x\in X$ with $\lda_x=+\infty$.
\end{defn}

In \cite{GL00,GL03}, this convergence was called {\em ball-affine convergence\/}
and it was proved that both the global and local versions
of this convergence coincide with the well-known slice convergence on the
space of \textit{convex} lsc functions. This
justifies the use in the present paper of the alternative name `slice' for this
convergence on the space of \textit{all\/} lsc functions.

The precise link between Wijsman convergence and slice convergence is described
in the theorem below whose proof is based on the following lemma:
\begin{lem}\label{permutx*}
Let $X$ be a normed space, $f,g:X\to\Rex$ and
$x^*\in X^*$. Then, $\gap(\graph g,\epi (f-x^*))>0$ if and
only if $\gap(\graph (g+x^*),\epi f)>0$.
\end{lem}
\pf
It is clearly sufficient to prove that the first condition implies
the second one. So, let $\eps>0$ such that $\gap(\graph g,\epi (f-x^*))
>\eps$, and then let $\delta\in (0,\eps/(1+\|x^*\|)$.
Pick $(x,\alpha)$ in $\epi f$ and $(y,\beta)$ in $\epi (g+x^*)$.
The lemma will be proved by showing that
$d((x,\alpha),(y,\beta))\geq \delta$.
The case $\|y-x\|\geq \delta$ being obvious, assume
$\|y-x\|< \delta<\eps$. Since
$(x, \alpha -\la x^*, x \ra)$ is in $\epi (f-x^*)$ and 
$(y,\beta-\la x^*, y \ra)$ is in $\graph g$, our assumption implies
that $\alpha -\la x^*, x \ra>\beta-\la x^*, y\ra +\eps$, hence
$\alpha -\beta>\eps+\la x^*, x-y\ra>\eps-\delta\|x^*\|> \delta$.
The proof is complete. 
\finpf

\begin{theorem}\label{W-versus-slice}
Let $X$ be a normed space and $\{f_n,f\}\subset\LSC(X)$.
The sequence $(f_n)_n$ is slice convergent to $f$ at $x$ with radius $\lda_x>0$
if and only if
every sequence $(f_n+x^*)_n$ with $x^*\in X^*$ is
Wijsman convergent to $f+x^*$ at $x$ with radius $\lda_x$.
\end{theorem}
\pf
It suffices to combine
Definition~\ref{variationalconv-slice},
Proposition~\ref{hit-miss-varcv} and Lemma~\ref{permutx*}.
\finpf

\section{Wijsman convergence and uniform infimum}
\label{WconvRobust}

The following analytic characterization of Wijsman convergence of functions
in terms of the lower limit of their infima on balls is an adaptation
of results in \cite{GL00,GL03}:

\begin{theorem} \label{anal-char} 
Let $X$ be a normed space and $\{f_n,f\}\subset\LSC(X)$.
Then, $(f_n)_n$ is Wijsman convergent to $f$ at $x$ with radius
$\lda_x\in (0,+\infty]$
if and only if for every $\lambda\in [0,\lda_x)$,
\smallskip\\
$ {\rm (W_a)} \left\lbrace 
    \begin{array}{rl} 
 {\rm (i)}  & \exists\, (x_n)_n\subset X: \;  
  x_n \to x, \; f_n(x_n) \to f(x);\smallskip\\ 
  {\rm (ii)} &   
    \dis r_{B_\lambda(x)}(f):=\sup_{\delta>0} \inf_{B_{\lambda+\delta}(x)}f 
    \leq \liminf_n \inf_{B_{\lambda}(x)}f_n. 
 \end{array} 
\right. 
$
\end{theorem} 

\pf
It suffices to show that `${\rm (W_g)(ii)}$ holds
for every $\lambda\in [0,\lda_x)$' is equivalent to
`${\rm (W_a)(ii)}$ holds for every $\lambda\in [0,\lda_x)$'.
To prove `$\imp$', fix $\lambda\in [0,\lda_x)$,
let $\alpha<\sup_{\delta>0} \inf_{B_{\lambda+\delta}(x)}f$,
and take $\gamma>0$ and $\delta\in (0,\gamma)$ such that
$$
\alpha+\gamma <\inf_{B_{\lambda+\delta}(x)} f.
$$
We have $\gap(B_{\lambda}(x)\times\alpha,\epi f) \geq \delta>0$,
because for
$z\in B_{\lambda}(x)$ and $(y,\beta) \in \epi f$ with
$\|y-z\|< \delta$, one has $\|y-x\| \leq \lambda +\delta$, so
$ \beta  \geq f(y) \geq \alpha+\gamma \geq  \alpha+\delta$,
whence $d((z,\alpha),(y,\beta))\ge \beta - \alpha  \geq  \delta$.
Applying ${\rm (W_g)\,(ii)}$, we derive that there exists $\delta>0$
such that
$\gap(B_{\lambda+\delta}(x)\times\alpha,\epi f_n) >0$ eventually,
so $f_n(z)>\alpha$ for every $z$ in $B_{\lambda+\delta}(x)$
and large $n$,
hence
$\liminf_n \inf_{B_{\lambda+\delta}(x)} f_n \geq \alpha.$
A fortiori, $\liminf_n \inf_{B_{\lambda}(x)} f_n \geq \alpha.$
The proof that ${\rm (W_a)\,(ii)}$ holds is complete.

\smallbreak
To prove `$\pmi$', fix $\lambda\in [0,\lda_x)$, let $\alpha\in\R$
such that $\gamma:=\gap(B_\lda(x)\times\alpha,\epi f) >0$, and
take $\delta >0$ such that
$$\gap(B_{\lda+\delta}(x)\times\alpha,\epi f)=\eps>0.$$
Let $z \in B_{\lda+\delta}(x)\cap \dom f$.
Since $(z,\alpha) \not \in \epi f$, one has $f(z) - \alpha>0$.
Therefore,
$$
f(z)-\alpha =d((z,\alpha),(z,f(z)))\geq
\gap(B_{\lda+\delta}(x)\times\alpha,\epi f)=\eps,
$$
so $f(z)\geq\alpha+\eps$ for every $z \in B_{\lambda+\delta}(x)$.
Hence, $\inf_{B_{\lambda+\delta}(x)}f>\alpha$.
Take $\deltb\in (0,\delta)$ such that $\ldab:=\lambda+\deltb<\lambda_x$.
Then,
$$
\sup_{\delta'>0}\inf_{B_{\ldab+\delta'}(x)}f
\ge \inf_{B_{\lambda+\delta}(x)}f>\alpha.
$$
\if{ Thus,
$$
\sup_{\delta>0}\inf_{B_{\lambda+\delta}(x)}f>\alpha.
$$
}\fi
Now, applying ${\rm (W_a)\,(ii)}$ with $\ldab$ instead of $\lda$, we find that
$\inf_{B_{\ldab}(x)} f_n=\inf_{B_{\lambda+\deltb}(x)} f_n>\alpha$
eventually, hence
$\epi f_n\cap(B_{\lambda+\deltb}(x)\times\alpha)=\emptyset$ 
eventually, from which we derive that, for any $\delta\in(0,\deltb)$,
$\gap(B_{\lambda+\delta}(x)\times\alpha,\epi f_n)\ge\deltb-\delta>0$
eventually, so that
$$
\sup_{\delta>0}
\liminf_n\gap(B_{\lambda+\delta}(x)\times\alpha,\epi f_n)>0,
$$
as required.
\finpf

The value $r_{B_\lambda(x)}(f):=\sup_{\delta>0} \inf_{B_{\lambda+\delta}(x)}f$
on the left of ${\rm (W_a)\,(ii)}$ cannot be replaced by the
usual infimum $\inf_{B_{\lambda}(x)}f$ when $f$ is only lsc
(see Example \ref{nogoodinf} below).
This is not fortuitous.
In problems involving non regular functions $f:X\to\Rex$ to be minimized
on a constrained set $S$, the value that naturally comes to the fore is
\begin{equation}\label{uniinf}
r_S(f) := \sup_{\delta>0}\inf_{B_\delta(S)}\, f
	= \lim_{\delta\to 0}\,\inf\{\, f(x)\tq d_S(x)\le \delta\,\}.
\end{equation}
The first explicit mention of (a variant of) this value dates back to
\cite{CL94}.
The value as written in \eqref{uniinf} was introduced and used in \cite{ACL99}.
Its importance was emphasized in \cite{Las01}, where the concept was generalized
and employed in various situations related to constrained minimization.
In the process, further properties and applications have been developed in
\cite{GL00,Jul03}. 
The notation $r_S(f)$ comes from \cite{ACL99}, slightly modifying the one in
\cite{CL94}. The name {\em uniform infimum of $f$ on $S$\/} for $r_S(f)$
was proposed in \cite{Las01}, arguing that this value
incorporates the behavior of $f$ on uniform neighborhoods of $S$.
Since then, this concept has been used in the textbooks \cite{BZ05,Pen13}
and in the survey \cite{Iof12} under different notations and names.
For example, in \cite{Pen13}, $r_S(f)$ is denoted
$\wedge_S(f)$ (more or less as in \cite{BZ05}) and is called
\textit{stabilized infimum};  
the usual infimum $\inf_S f$ is declared \textit{robust} when it is equal
to $r_S(f)$; a point $\xb$ achieving this value, i.e.\ $f(\xb)=r_S(f)$,
is called a \textit{robust minimizer} (more or less as in \cite{Iof12}).

In general, $r_S(f)<\inf_S f$ for arbitrary lsc $f$:
additional conditions (so-called qualification conditions) are required
to have the equality $r_S(f)=\inf_S f$.
\medbreak\noindent
\begin{exe}\label{nogoodlsc}
A lsc $f$ with $r_{B_\lambda(0)}(f)<\inf_{B_\lambda(0)}f$
for arbitrary small $\lambda$
(see also \cite[Example 2]{Las01}, \cite[Exemple 3.7]{Jul03}).
Let $X$ be the Hilbert space $\ell^2(\N)$ and let $(e_i)_{i\in\N}$
be its canonical basis. Define
$$
f(x):=\left\{
	\begin{array}{ll}
	-1/n & \mbox{if }~ x=e_i/n + e_1/in, ~~n=1,2,\ldots, ~~i=1,2,\ldots\\
	0  & \mbox{otherwise}.
	\end{array}
      \right .
$$
Then $f$ is lsc at every point. For $\lambda=1/n$, $n\in\N$,
one has
$$
r_{B_\lambda(0)}(f)=-1/n, \quad \inf_{B_\lambda(0)}f=-1/(n+1).
$$
So, for every $\lambda>0$ there exists $\lambda_n\in (0,\lambda)$
such that
$r_{B_{\lambda_n}(0)}(f)<\inf_{B_{\lambda_n}(0)}f.$
\end{exe}

\begin{exe}
Sufficient conditions for $r_S(f)=\inf_S f$ to hold
(see \cite[Proposition 3.2]{ACL99}):

1. $S=X$,

2. $f$ is uniformly continuous on a uniform neighborhood of $S$,

3. $f$ is lsc on a neighborhood of $S$
and $S$ is compact, or $f$ is inf-compact and $S$ is closed,

4. $X=\R_+(\dom f - S)$, $f$ is convex lsc and $S$ is closed and convex.
\end{exe}

Since $r_S(f)$ is the natural value bound to the constrained minimization
problem\smallskip\\
\centerline{($\mathcal{P}$)  Minimize $f(x)$ subject to $x\in S$}
\smallbreak\noindent
it is expected that $r_S(f)$ can be obtained as the limiting value of the unconstrained penalized problems associated with ($\mathcal{P}$).
This is indeed the case.
The following proposition was established in \cite[Proposition 3.16]{Jul03}
gathering earlier observations
(see also \cite[Proposition 1.130]{Pen13}).
For the sake of completeness, we reproduce the proof.  

\begin{prop}\label{dec-inf}
Let $X$ be a normed space, $f:X\to\Rex$ bounded from below and
$S\subset X$ such that $S\cap \dom f\ne\emptyset$.
Then, for any $p>0$,
\begin{equation}\label{dec-inf1}
r_S(f)=
\lim_{n\to\infty}\,\inf\{f(x)+n d^p_S(x)\tq x\in X\}.
\end{equation}
\end{prop}

\pf
Let $\eps>0$ and $\eta>0$. Choose $\delta >0$ such that
$\eta\delta^p<\eps$. Then,
$$
\inf_X \,(f+\eta d^p_S)\le \inf_{B_\delta(S)}f +\eta\delta^p
\le \sup_{\delta>0} \inf_{B_\delta(S)}f +\eps,
$$
showing that the first member of (\ref{dec-inf1}) is not smaller than
the second one.

Now, let $\gamma<\sup_{\delta>0} \inf_{B_\delta(S)}f$. Take
$\delta>0$ such that  $\gamma<\inf_{B_\delta(S)}f$ and choose $\eta>0$
such that
$\inf_{B_\delta(S)}f\le \eta\delta^p+\inf_Xf$
(this is possible since both $\inf_{B_\delta(S)}f$ and $\inf_Xf$
are finite).
We claim that
$$
\inf_{B_\delta(S)}f\le f(x)+\eta d^p_S(x),\quad \forall x\in X.
$$
This is clear if $x$ belongs to $B_\delta(S)$;
otherwise $d_S(x)\ge \delta$, hence, due to our choice of $\eta$,
$\inf_{B_\delta(S)}f\le \inf_Xf + \eta\delta^p\le f(x)+\eta d^p_S(x)$. 
It follows that
$$
\gamma <\sup_{\eta>0}\inf_X \,(f+\eta d^p_S),
$$
showing that the first member of (\ref{dec-inf1}) is not greater than
the second one.
\finpf

The next two propositions provide useful examples of
Wijsman convergent sequences.

\begin{prop}\label{carac-W}
Let $X$ be a normed space, $f\in LSC(X)$ and
$S\subset X$ a closed subset such that $S\cap \dom f\ne\emptyset$.
Let $f_n:=f+nd_S^p$ with $p>0$.
The following are equivalent:

{\rm (a)} The sequence $(f_n)_n$ is Wijsman convergent to $f_S$
at every point $x\in X$,

{\rm (b)} $r_{B_\lambda(x)}(f_S)\le r_S(f_{B_\lambda(x)})$
for every $x\in X$ and $\lambda>0$ small enough.
\end{prop}
\pf
We have to show that the assertions ${\rm (W_a)}$ in Theorem \ref{anal-char}
hold if and only if (b) holds.
Assertion ${\rm (W_a)(i)}$ is always satisfied at
every point $x\in X$ by the constant sequence $x_n:=x$. Indeed,
$f_n(x)=f_S(x)$ for $x\in S$ and $\lim_n f_n(x)=+\infty$
for $x\not\in S$ since $S$ is closed,
so $\lim_n f_n(x)=f_S(x)$ for $x\not\in S$.
On the other hand, 
${\rm (W_a)(ii)}$ asserts that for $\lambda>0$ small enough,
$$
r_{B_\lambda(x)}(f_S)
\le\liminf_{n\to\infty}\,\inf \{f_n(z)\tq z\in {B_\lambda(x)}\}
=\lim_{n\to\infty}\,\inf \{f_{B_\lambda(x)}(z)+nd_S(z)\tq z\in X\}.
$$
Now by lower semicontinuity, $f_{B_\lambda(x)}$ is bounded from below
for $\lambda>0$ small enough, so  according to Proposition \ref{dec-inf},
the expression on the right hand side is equal to $r_S(f_{B_\lambda(x)})$
for $\lambda>0$ small enough.
Hence, ${\rm (W_a)(ii)}$ is equivalent to (b).
\finpf

\begin{prop}\label{infconvWijsman}
Let $X$ be a normed space and let $f\in LSC(X)$ be proper and bounded
from below.
Let $f_n:=f\triangledown n\|.\|$. 
Then, each $f_n$ is Lipschitz continuous on $X$
and the sequence $(f_n)_n$ is Wijsman convergent to $f$ with
$r_{B_\lambda(x)}(f)=\lim_n \inf_{B_\lambda(x)}f_n$ for every $x\in X$ and
$\lambda \ge 0$.
\end{prop}
\pf
We have $f_n(x)= \inf_{y\in X}(f(x+y)+n\|y\|)$.
We first observe that $f_n(x)$ is finite for every $x\in X$ and $n\in \N$.
Indeed, $f$ being proper, there is $y\in X$ such that
$f(x+y)$ is finite, so $f_n(x)<+\infty$, and $f$ being bounded from below,
$-\infty<\inf_X f\le f_n(x)$.
Otherwise, for any $x, v\in X$ and $n\in \N$, it holds
\begin{equation}\label{Lip}
f_n(x)\le \inf_{y\in X}(f(x+y)+n\|y-v\|)+ n\|v\|
=f_n(x+v)+n\|v\|,
\end{equation}
so $|f_n(x)-f_n(x+v)|\le n\|v\|$, proving that each $f_n$ is $n$-Lipschitz.
\smallbreak
We now prove the second statement. For any $\lambda\ge 0$, one has
$$
\inf_{B_\lambda(x)}f_n
=\inf_{x'\in B_\lambda(x)}\inf_{z\in X}(f(z)+n\|z-x'\|)
=\inf_{z\in X}(f(z)+nd_{B_\lambda(x)}(z)).
$$
So, taking the limit as $n\to +\infty$ on both sides,
we see from Proposition \ref{dec-inf} that
\begin{equation}\label{toto}
\lim_{n\to\infty} \inf_{B_\lambda(x)}f_n=r_{B_\lambda(x)}(f).
\end{equation}
For $\lambda= 0$,
\eqref{toto} gives $\lim_{n\to\infty} f_n(x)=r_{\{x\}}f$,
where $r_{\{x\}}f=\sup_{\delta>0}\inf_{B_\delta(x)}f=f(x)$
because $f$ is lsc, so ${\rm (W_a)(i)}$
is satisfied at any point $x\in X$ by the constant sequence $x_n:=x$.
On the other hand, \eqref{toto} clearly implies ${\rm (W_a)(ii)}$.
This proves that the sequence $(f_n)_n$ is Wijsman convergent to $f$.
\finpf

\begin{exe}\label{nogoodinf}
Let $f$ be the lsc function defined in Example \ref{nogoodlsc}
and let $(f_n)_n$ be the sequence given in Proposition \ref{infconvWijsman}.
This sequence is Wijsman convergent to $f$ with
$r_{B_\lambda(x)}(f)=\lim_n \inf_{B_\lambda(x)}f_n$ for every $x\in X$ and
$\lambda \ge 0$.
Since for every $\lda_x>0$ there exists $\lambda\in [0,\lda_x)$
such that $r_{B_\lambda(x)}(f)< \inf_{B_\lambda(x)}f$, we infer that
for every $\lda_x>0$ there exists $\lambda\in [0,\lda_x)$
such that $\liminf_n \inf_{B_\lambda(x)}f_n< \inf_{B_\lambda(x)}f$.
This shows that we cannot replace the value $r_{B_\lambda(x)}(f)$
by $\inf_{B_\lambda(x)}f$ in ${\rm (W_a)(ii)}$.
\end{exe}
\section{Stability of slopes with respect to Wijsman convergence}
\label{slopestability}
From now on, $X$ denotes a Banach space and $f_n,f:X\to \Rex$ denote
\lsc\ functions.
The {\em slope\/} of $f:X\to\Rex$ at $x\in \dom f$, introduced in
\cite{DMT80}, is defined by
$$
|\nabla f|(x):= \limsup_{\substack{y \to x \\ y\not =x}}
\frac{(f(x)-f(y))^+} {\|x-y\|},
$$
where $\alpha^+=\max(0,\alpha)$ for $\alpha\in\Rex$.

\begin{theorem}\label{slope-stab}
If $(f_n)_n$ is Wijsman convergent to $f$ at $x\in \dom f$, then
there is a sequence $(x_n)_n\subset X$ such that
$x_n\to x$, $f_n(x_n)\to f(x)$, and
\begin{equation*}\label{W-slope}
\limsup_n\,|\nabla f_n|(x_n)\le|\nabla f|(x).
\end{equation*}
\end{theorem}
\pf
The proof is exactly the same as the one of \cite[Theorem 3.1]{GL03},
where slice convergence was considered instead of Wijsman convergence,
but we reproduce it for the reader's convenience.
Let $\sigma:=|\nabla f|(x)$, which we may suppose to be finite.
Let $(x_n)_n\subset X$ be a sequence verifying ${\rm (W_a)\,(i)}$
and let $\ldx> 0$ be such that ${\rm (W_a)\,(ii)}$ holds for
$\lda\in[0,\ldx)$.
We claim that
for every $\eps\in (0,\ldx)$ there exists $N_\eps\in\N$ such
that for each $n\ge N_\eps$ there exists
$\bx_n\in X$ verifying
\begin{equation}\label{claim}
\left\lbrace \begin{array}{l} 
\|\bx_n-x_n\|<\eps, \quad
      |f_n(\bx_n)- f_n(x_n)|\leq (\sigma+2\eps)\eps,
\smallbreak\\ 
|\nabla f_n| (\bx_n)\leq \sigma + 3\eps.
\end{array}\right.
\end{equation}
Indeed, it follows from the definition of the slope of $f$ at $x$ that
there exists $\ldx'\in (0,\eps)$ such that
for all $\lambda\in (0,\ldx')$,
$$f(x)<\inf_{B_\lambda(x)} f + (\sigma+\eps)\lambda,$$
hence, 
for all $\lambda\in (0,\ldx')$,
$$
f(x)<\sup_{\delta> 0}\inf_{B_{\lambda+\delta}(x)}f+(\sigma+\eps)\lambda.
$$
Fix $\lambda\in (0,\ldx')$. Combining the previous inequality with
${\rm (W_a)\,(ii)}$ we get
$$
f(x)<\liminf_n \inf_{B_{\lambda}(x)}f_n +(\sigma+\eps)\lambda,
$$
while, according to ${\rm (W_a)\,(i)}$, for $n$ large enough one has
$$
f_n(x_n) < f(x) + \lambda\eps,
$$
$$
B_\mu (x_n) \subset B_\lambda(x),\quad
\mbox{with } \mu:=\frac{\sigma+2\eps}{\sigma+3\eps}\lambda.
$$
We thus derive that for each $n$ large enough one has
$$
f_n(x_n)< \inf_{B_\mu(x_n)} f_n + (\sigma+2\eps)\lambda
         =\inf_{B_\mu(x_n)} f_n + (\sigma+3\eps)\mu.
$$
Then, applying Ekeland's variational principle, we get a point
$\bx_n \in X$ such that
$$\left\{
\begin{array}{l}
\|\bx_n-x_n\|<\mu,
\quad f_n(\bx_n)\leq f_n(x_n)\leq f_n(\bx_n)+(\sigma+2\eps)\lambda,
\smallskip\\
f_n(\bx_n)\le f_n(x)+(\sigma+3\eps) \|x-\bx_n\|,
           \quad\forall x\in B_\mu(x_n),
\end{array}
\right.
$$
which implies that $|\nabla f_n|(\bx_n) \leq \sigma+3\eps$.
The proof of the claim is therefore complete.

Next, for each $\eps=1/k<\ldx$, choose an integer $N_k$ and
a sequence $(\bx_{n,k})_n$ in $X$ verifying (\ref{claim})
for each $n\ge N_k$. Without loss of generality, we may assume
that $N_{k+1}> N_k$. The desired sequence is then given by
$\bx_n:=\bx_{n,k}$ for $n\in \N$ and $k$ such that $N_k\le n< N_{k+1}$.
\finpf

\begin{cor}\label{cor1}
Assume $|\nabla f|(\xb)=0$.
If $(f_n)_n$ is Wijsman convergent to $f$ at $\xb$, then
there is a sequence $(x_n)_n$ such that $x_n\to \xb$,  $f_n(x_n)\to f(\xb)$
and $|\nabla f_n|(x_n)\to 0$.
\end{cor}

\begin{cor}\label{cor2}
Assume $\inf_X f>-\infty$. 
If $(f_n)_n$ is Wijsman convergent to $f$, then
there is a sequence $(x_n)_n$ such that $f_n(x_n)\to \inf_X f$
and $|\nabla f_n|(x_n)\to 0$.
\end{cor}
\pf
For each positive integer $n$, let $y_n\in X$ such that
$f(y_n)\le \inf_X f+1/n^2$. Apply Ekeland's variational principle
to get $z_n\in X$ such that $f(z_n)\le f(y_n)\le \inf_X f+1/n^2$ and
$f(z_n)\le f(y)+(1/n)\|y-z_n\|$ for every $y\in X$.
This implies that $|\nabla f|(z_n)\le 1/n$. Now, by Theorem \ref{slope-stab},
we can construct a sequence $(x_n)_n$ such that for each $n$,
$|f_n(x_n)-f(z_n)|\le 1/n$ and $|\nabla f_n|(x_n)\le |\nabla f|(z_n)\le 1/n$.
This sequence has the required properties.
\finpf

\section{Trustworthiness and stability of subdifferentials}

In the sequel, we call \textit{subdifferential} any operator $\del$ that
associates a set $\del f(x)\subset X^*$ to any triplet $(X,f,x)$,
where $X$ is Banach space,
$f\in\LSC(X)$ and $x \in X$, in such a way that
the following properties are satisfied:\smallbreak
(A1) If $f$ is convex near $x$, then
$\del f(x) =\{ x^*\in X^* \tq
       \langle x^*, y-x \rangle + f(x) \le f(y), ~~ \forall y\in X \}$;
\smallbreak
(A2) If $F(x,y) = f(x) + g(y)$, then 
$\del F(x,y)\subset \del f(x)\times\del g(y)$;\smallbreak
(A3) For any $f$ and $x^*\in X^*$, $\del (f+x^*)(x)=\del f(x)+x^*$.
\smallbreak\noindent
There are many other basic properties shared by all interesting
subdifferentials (see \cite[Definition 2.1]{Iof12}). But in what follows
we need only the three properties above.
\medbreak
We write
$\del f:=\{(x,x^*)\in X\times X^*\tq x^*\in \del f(x)\}$
for the graph of $\del f$.
As in \cite[Definition 2.12]{Iof12}, we say that
a subdifferential $\del$ is \textit{trustworthy} on a space $X$,
or that $X$ is a \textit{trustworthy} space for $\del$,
if the following rule holds:
\smallbreak\noindent
(R1) \textit{Fuzzy minimization rule.}
For any $f\in \LSC(X)$ and $\varphi$ convex Lipschitz,
if $f+\varphi$ admits a finite local minimum at $z$, then
there are sequences $((x_n,x_n^*))_n\subset \partial f$ and
$((y_n,y_n^*))_n\subset \partial \varphi$ such that
$x_n\to z$, $y_n\to z$, $f(x_n)\to f(z)$ and
$x_n^*+y_n^*\to 0$.
\medbreak\noindent
\begin{exe}
Main trustworthy spaces
(see \cite{Las01,BZ05,Iof12,Pen13} and the references therein):

1. $X$ is a Hilbert space, $\del$ is the proximal subdifferential;

2. $X$ is an Asplund space, $\del$ is
the Fr\'echet or the limiting Fr\'echet subdifferential;

3. $X$ is a separable Banach space, $\del$ is
the Hadamard subdifferential;

4. $X$ is any Banach space, $\del$ is the
subdifferential of Clarke, of Michel-Penot, or of Ioffe.
\end{exe}
\medbreak
The rule (R1) for trustworthiness expresses
a subdifferential necessary condition for a point $z$
to be a local minimizer of the penalized function $f+\varphi$
where $\varphi$ is convex Lipschitz.
In fact, trustworthiness can be characterized by various properties
related to such penalized functions.
\if{In fact, trustworthiness characterizes subdifferentials that
behave well when facing convex Lipschitz penalizations.}\fi
For example:

\medbreak\noindent
(P1) 
\textit{Necessary condition for an approximate local minimizer.}
For any $f\in \LSC(X)$, $\varphi$ convex Lipschitz,
$\lambda>0$ and $\sigma>0$,
if\smallskip\\
\centerline{
$(f+\varphi)(z) < \inf_{B_\lambda(z)} (f+\varphi) + \lambda\sigma,$}
\smallbreak\noindent
then
there exist $(\xb,\xb^*) \in \del f$ and $(\yb,\yb^*)\in \del \varphi$
such that
$\|\xb-z\| < \lambda$,  $\|\yb-z\|<\ld$,
$|f(\xb)+\varphi(\yb)-(f(z)+\varphi(z)| < \lambda\sigma$
and $\|\xb^*+\yb^*\| < \sigma$.
 
\bigbreak\noindent
(P2) \textit{Slope control.} 
For any $f\in \LSC(X)$, $\varphi$ convex Lipschitz and $z\in\dom f$,
there are sequences $((x_n,x_n^*))_n\subset \partial f$ and
$((y_n,y_n^*))_n\subset \partial \varphi$ such that\smallskip\\
\centerline{
$x_n\to_f z$, $y_n\to z$ and
$\limsup_n\|x_n^*+y_n^*\|\le |\nabla (f+\varphi)|(z)$.}


\bigbreak\noindent
(P3) \textit{Fr\'echet subdifferential control.} 
For any $f\in \LSC(X)$, $\varphi$ convex Lipschitz and
$(z,z^*)\in \del_F (f+\varphi)$,
there are sequences $((x_n,x_n^*))_n\subset \partial f$ and
$((y_n,y_n^*))_n\subset \partial \varphi$ such that\smallskip\\
\centerline{
$x_n\to_f z$, $y_n\to z$ and $x_n^*+y_n^*\to z^*$.}

\bigbreak
We recall that the {\em Fr\'echet subdifferential\/} of $f$ at $\xb$ is
given by 
\begin{equation}\label{FreDef}
\del_F f(\xb):=\{\,\xbs \in X^* \tq 
   \liminf_{\substack{x\to\xb \\ x\not = \xb}}\,
\frac {f(x)-f(\xb)-\langle\xbs,x-\xb\rangle}{\|x-\xb\|}\geq 0\,\},
\end{equation}
which, as observed in \cite[Proposition 4.1]{Jul03},
can be conveniently rewritten as
\begin{equation}\label{FreDef2}
\del_F f(\xb)=\{\,\xbs\in X^*\tq |\nabla(f-\xbs)|(\xb)=0\,\}.
\end{equation}

\medbreak
Property (P1) was considered for the first time in \cite{GL00}
(in the special case where $\varphi$ is a continuous linear form)
and in \cite{MW02} for the Fr\'echet subdifferential and $\varphi=0$.
The following proposition was established in {\cite[Th\'eor\`eme 4.2]{Jul03}}.
We briefly recall the proof for the sake of completeness.

\begin{prop}\label{equiv-OP}
Rule (R1) and Properties (P1)--(P3) are equivalent.
\end{prop}
\pf
(R1) $\imp$ (P1) with $\varphi$ linear continuous
was already observed in \cite[Theorem 3.2]{GL00}.
Let $f\in \LSC(X)$, $\varphi$ convex Lipschitz,
$\lambda>0$ and $\sigma>0$ such that
$(f+\varphi)(z) < \inf_{B_\lambda(z)} (f+\varphi) + \lambda\sigma$.
Apply Ekeland's variational principle to $g := f + \varphi +
\delta_{B_\lambda(z)}$ with $0 < \sigma' < \sigma$ such that
$g(z) < \inf_X g + \lambda\sigma'$ and with $\lambda'$ such that
$0 < \lambda' < \lambda$.
We obtain a point $\zb \in B_{\lambda'}(z)$, with
$|(f+\varphi)(\zb) - (f+\varphi)(z)|\leq \lambda\sigma'$, that
is a local minimizer of the function
$f+\varphi+\sigma'\|\;.-\zb\|$.

Now we apply (R1) to $f$ and  
$\psi:=\varphi+\sigma'\|\;.-\zb\|$ at point $\zb$
to get $(\xb,\xbs) \in \del f$ and $(\yb,\ybs) \in \del \psi$
such that $\|\xb-\zb\| <\lambda- \lambda'$,
$|f(\xb)-f(\zb)|< \lambda(\sigma-\sigma')/2$,
$\varphi(\yb)-\varphi(z)|<\lambda(\sigma-\sigma')/2$ and
$\|\xbs+\ybs\| < \sigma-\sigma'$.
Combining the above inequalities, we infer that
$\|\xb-z\|< \lambda$ and $|f(\xb)+\varphi(\yb)-f(z)-\varphi(z)|< \lambda\sigma$.
On the other hand, as $\ybs \in \del \psi(\yb)$,
using (A1) and standard calculus rules of convex analysis,
we derive that $\ybs=\ybs_0+\sigma'\xi^*$
where $\ybs_0\in \del\varphi(\yb)$ and $\|\xi^*\|\le 1$, so
$$\|\xbs+\ybs_0\| \leq \|\xbs+\ybs\| + \|\ybs_0-\ybs\\| < \sigma.$$
Thus, $(\xb,\xbs) \in \del f$ and $(\yb,\ybs_0)\in \varphi$
satisfy the required inequalities of (P1).
\smallbreak
(P1) $\imp$ (P2). Let $|\nabla (f+\varphi)|(z)<\sigma$.
We have to show that for any $\eps>0$, there exist
$(x,x^*)\in \del f$ and $(y,y^*)\in\del\varphi$ such that
$\|x-z\|< \eps$, $\|y-z\|< \eps$, $|f(x)-f(z)|<\eps$
and $\|x^*+y^*\|<\sigma$.
Fix $\sigma'$ such that $|\nabla (f+\varphi)|(z)<\sigma'<\sigma$.
From the definition of the strong slope, we derive that there is
$\ld'>0$ such that $(f+\varphi)(x)-(f+\varphi)(y)<\sigma'\|x-y\|$
for every $y\in B_{\ld'}(x)$, which implies that
for all $\ld\in]0,\ld']$, one has
$(f+\varphi)(x)<(f+\varphi)(y)+\sigma'\ld$ for every $y\in B_{\ld}(x)$.
Hence, for every $\ld\in]0,\ld']$, it holds
$$
(f+\varphi)(x)\le \inf_{B_{\ld}(x)} (f+\varphi)+\sigma'\ld
      <\inf_{B_{\ld}(x)} (f+\varphi)+\sigma\ld.
$$
Let $\eps>0$. Apply (P1) with $\ld=\min\{\ld',\eps,\eps/\sigma\}$. We obtain
$(x,x^*)\in \del f$ and $(y,y^*)\in\del\varphi$ such that
$\|x-z\|<\ld\le\eps$, $\|y-z\|<\ld\le\eps$,
$|f(x)+\varphi(y)-f(z)-\varphi(z)|<\ld\sigma\le\eps$
and $\|x^*+y^*\|<\sigma$.
Since $\varphi$ is continuous, we can manage so that the penultimate
inequality induces $|f(x)-f(z)|<\eps$. The proof is complete.
\smallbreak
(P2) $\imp$ (P3). Let $z^*\in \del_F (f+\varphi)(z)$. By \eqref{FreDef2},
this amounts to $|\nabla(f+\varphi-z^*)|(z)=0$.
Applying (P2) we get sequences $((x_n,x_n^*))_n\subset
\partial (f-z^*)$ and $((y_n,y_n^*))_n\subset \partial \varphi$ 
with
$x_n\to z$, $y_n\to z$, $(f-z^*)(x_n)\to (f-z^*)(z)$ and $x_n^*+y_n^*\to 0$.
From (A3), we derive that $\xb_n^*:=x_n^*+z^*\in \del f(x_n)$.
Then, the sequences $((x_n,\xb_n^*))_n$ and $((y_n,y_n^*))_n$
satisfy the required properties.
\smallbreak
(P3) $\imp$ (R1) is obvious.
\finpf

The following theorem asserts that Property (P2) is stable with respect
to Wijsman perturbations, the next one that Property (P3) is stable with
respect to slice perturbations.

\begin{theorem}\label{subdifWstable}
Assume (P2) holds on the space $X$. Let $(f_n)_n$ be a sequence of lsc functions
and $(\varphi_n)_n$ a sequence of convex Lipschitz functions such that
$(f_n+\varphi_n)_n$ is Wijsman convergent to a function $f$ at $\zb$.
Then, there are elements $(x_n,x_n^*)\in \partial f_n$ and
$(y_n,y_n^*)\in \partial \varphi_n$
such that $x_n\to \zb$,  $y_n\to \zb$,
$f_n(x_n)+\varphi_n(y_n)\to f(\zb)$
and $\limsup_n \|x^*_n+y_n^*\|\le |\nabla f|(\zb)$.
\end{theorem}

\pf
By Theorem \ref{slope-stab}, there is a sequence $(z_n)_n$ such that
$z_n\to \zb$,  $(f_n+\varphi_n)(z_n)\to f(\zb)$ and
$\limsup_n |\nabla (f_n+\varphi_n)|(z_n)\le |\nabla f|(\zb)$.
By (P2) applied to $f_n$, $\varphi_n$ and $z_n$,
there are sequences $(x_n)_n$ and $(y_n)_n$ such that,
for every $n\in \N$, $\|x_n-z_n\|\le 1/n$, $\|y_n-z_n\|\le 1/n$
$|f_n(x_n)-f_n(z_n)|\le 1/n$ and 
$$
d\left(0, \del f_n(x_n)+\del \varphi_n(y_n)\right)
< |\nabla (f_n+\varphi_n)|(z_n)+1/n.
$$
So, for every $n\in \N$, there are elements $x_n^*\in \del f_n(x_n)$
and $y_n^*\in \del \varphi_n(y_n)$ such that
$$\|x_n^*+y_n^*\|<|\nabla (f_n+\varphi_n)|(z_n)+1/n.$$
Since $z_n\to \zb$,  $f_n(z_n)\to f(\zb)$ and
$\limsup_n |\nabla (f_n+\varphi_n)|(z_n)\le |\nabla f|(\zb)$,
the sequences $((x_n,x_n^*))_n$ and $((y_n,y_n^*))_n$
satisfy the required properties.
\finpf

\begin{theorem}\label{jetSlicestable}
Assume (P2) holds on the space $X$.
Let $(f_n)_n$ be a sequence of lsc functions
and $(\varphi_n)_n$ a sequence of convex Lipschitz functions such that
$(f_n+\varphi_n)_n$ is slice convergent to a function $f$. Then,
for each $(\zb,\zb^*)\in \del_F$,
there are elements $(x_n,x_n^*)\in \partial f_n$ and
$(y_n,y_n^*)\in \partial \varphi_n$
such that $x_n\to z$,  $y_n\to z$,
$f_n(x_n)+\varphi_n(y_n)\to f(z)$ and $x^*_n+y_n^*\to z^*$.
\end{theorem}

\pf 
Let $(\zb,\zb^*)\in \del_F$. So $|\nabla(f-\zb^*)|(z)=0$.
By Theorem \ref{W-versus-slice}, the sequence $((f_n-\zb^*+\varphi_n)_n$ is
Wijsman convergent to $f-\zb^*$.
Applying Theorem \ref{subdifWstable} and (A3), we find
sequences $((x_n,x_n^*))_n$ and $((y_n,y_n^*))_n$ such that
$x_n\to \zb$,  $y_n\to \zb$,
$f_n(x_n)+\varphi_n(y_n)\to f(\zb)$
and $\|x^*_n-\zb^*+y_n^*\|\to 0$. The conclusion follows.
\finpf

Variants of Theorem \ref{subdifWstable} and Theorem \ref{jetSlicestable}
are considered in \cite[Theorem 4.1]{GL00} and \cite[Theorem 5.1]{GL03}
and several applications are given.

\section{Subdifferential Sum Rules}
\label{subdiffsum}
A family of lsc functions $\{f_1,\ldots,f_k\}$ is said to be
{\em decouplable\/} at $\xb\in X$ (\cite[D\'efinition 3.5]{Jul03})
provided there is $\lda_x>0$ such that for any $\lda\in [0,\lda_x)$,
\begin{equation}\label{CD}
\hspace*{-7pt}
\sup_{\delta>0}\inf_{B_{\lda+\delta}(\xb)}\sum f_i\le
\sup_{\delta>0}\,\inf\left\{\,\sum f_i(x_i)\tq x_i\in B_{\lda}(\xb),
\|x_i-x_j\|\le\delta\,\right\}.
\end{equation}
It is said to be \textit{$X^*$-decouplable} at $\xb$ if
$\{f_1,\ldots,f_k,x^*\}$ is decouplable at $\xb$ for every $x^*\in X^*$.
\begin{exe}
Sufficient conditions for a family $\{f_1,\ldots,f_k\}$ to be
decouplable at $\xb$:

1.  All but at most one of the functions are uniformly continuous near $\xb$
\cite{Las01}.

2. At least one of the functions has compact lower level sets near $\xb$
\cite{Las01}.

3. The function $\sum f_i$ achieves a local uniform (decoupled, robust)
minimum at $\xb$ \cite{Las01,BZ05,Iof12}.

4. $k=2$ and the inf-convolution
$f_{B_{\lda}(\xb)}\triangledown g^{-}_{B_{\lda}(\xb)}$
(with $f=f_1$, $g=f_2$) is lsc at $0$ \cite{Jul03}.
\end{exe}

Let $F:(x_1,\ldots,x_k)\in X^k\mapsto \sum f_i(x_i)$
be the \textit{decoupled sum},
$B^k_{\lda}(\xb):=B_{\lda}(\xb)\times\cdots\times B_{\lda}(\xb)$
be the $\lda$-ball of center $(\xb,\ldots,\xb)$ in $X^k$ with the max-morm,
and let $\Delta:=\{(x,\ldots,x)\in X^k\tq x\in X\}$ be the diagonal of $X^k$. 
The two expressions in the decoupling's condition \eqref{CD} can be written as
\smallbreak
$\dis
\sup_{\delta>0}\inf_{B_{\lda+\delta}(\xb)}\sum f_i
=r_{B_{\lda}(\xb)}\left(\sum f_i\right)=r_{B^k_{\lda}(\xb)}(F_\Delta),
$

\begin{multline*}
\sup_{\delta>0}\,\inf\left\{\,\sum f_i(x_i)\tq x_i\in B_{\lda}(\xb),
\|x_i-x_j\|\le\delta\,\right\}=\\
\sup_{\delta>0}\,\inf\left\{\,\sum {f_i}_{B_{\lda}(\xb)}(x_i)\tq
d_\Delta(x_1,\ldots,x_k)\le \delta\,\right\}
=r_\Delta(F_{B^k_{\lda}(\xb)}).
\end{multline*}
\noindent
So the decoupling's condition \eqref{CD} amounts to:

\begin{equation}\label{CDbis}
r_{B^k_{\lda}(\xb)}(F_\Delta)\le r_\Delta(F_{B^k_{\lda}(\xb)}).\tag{\ref{CD}b}
\end{equation}

\begin{prop}\label{decou-Wijs}
The family $\{f_1,\ldots,f_k\}$ is decouplable at $\xb\in X$ if and only if the
sequence $F_n:= F+nd_\Delta$ is Wijsman convergent to
$F_\Delta$ at $(\xb,\ldots,\xb)$.
\end{prop}

\pf
According to the above, the family
$\{f_1,\ldots,f_k\}$ is decouplable at $\xb\in X$ if and only if
\eqref{CDbis} holds for any $\lda>0$ small enough,
and according to Proposition \ref{carac-W}, this latter condition
means that the sequence $F_n:= F+nd_\Delta$ is Wijsman convergent to
$F_\Delta$ at $(\xb,\ldots,\xb)$.
\finpf

We consider generalizations of (P2) and (P3) to
decouplable families of functions:

\medbreak\noindent
(R2) \textit{Slope control.} 
Let $\{f_1,\ldots,f_k\}\subset \LSC(X)$ be decouplable at
$\xb\in X$.
If $|\nabla (\sum f_i)|(\xb)<\infty$,
then there are sequences $((x_{i,n},x_{i,n}^*))_n\subset\del f_i$,
$i=1,\ldots,k$, such that \smallskip\\
\centerline{$x_{i,n}\to_{f_i} \xb$,
$\limsup_n\|\sum x^*_{i,n}\|\le
\left |\nabla \left(\sum f_i\right)\right |(\xb)$
and $\diam (x_{1,n},\ldots,x_{k,n})\,\|x^*_{i,n}\|\to 0$.}


\bigbreak\noindent
(R3) \textit{Fr\'echet subdifferential control.} 
Let $\{f_1,\ldots,f_k\}\subset \LSC(X)$ be $X^*$-decouplable at
$\xb\in X$.
For any $\xb^*\in \del_F \left(\sum f_i\right)(\xb)$,
there are sequences $((x_{i,n},x_{i,n}^*))_n\subset\del f_i$,
$i=1,\ldots,k$, such that\smallskip\\
\centerline{
$x_{i,n}\to_{f_i} \xb$, $\sum x^*_{i,n}\to \xb^*$
and $\diam (x_{1,n},\ldots,x_{k,n})\,\|x^*_{i,n}\|\to 0$.}

\begin{theorem}\label{subsum}
Let $\mathcal{X}$ be a class of Banach spaces which contains Cartesian
products of its elements.
If (R1) holds on every space in $\mathcal{X}$,
then so do (R2) and (R3).
\end{theorem}

\pf
Let $X$ be a space in the class $\mathcal{X}$. We show with some details
that (R2) holds on $X$.
By Proposition \ref{decou-Wijs}, the sequence of functions $F+nd_\Delta$
defined on $X^k$ is Wijsman convergent to the function
$F_\Delta$ at $(\xb,\ldots,\xb)$.
Since (R1) holds on $X^k$, we may apply Theorem \ref{subdifWstable}
with $f_n:=F$, $\varphi_n:=nd_\Delta$, $f:=F_\Delta$ and $\zb:=(\xb,\ldots,\xb)$.
This produces  elements $(\hat{x}_n, \hat{x}^*_n)\in \del F$
and $(\hat{y}_n, \hat{y}^*_n)\in \del (nd_\Delta)$
such that 
\smallbreak
(a) $\hat{x}_n\to (\xb,\ldots,\xb)$, $\hat{y}_n\to (\xb,\ldots,\xb)$,
\smallbreak
(b) $F(\hat{x}_n)+nd_\Delta(\hat{y}_n)
\to F_\Delta(\xb,\ldots,\xb)=\sum f_i(\xb)$,
\smallbreak
(c) $\limsup_n \|\hat{x}^*_n+\hat{y}^*_n\|
\le |\nabla F_\Delta|(\xb,\ldots,\xb)
=\left |\nabla \left(\sum f_i\right)\right |(\xb)$.

\medbreak\noindent
The point $\hat{x}_n$ can be written as
$\hat{x}_n=(x_{1,n},\ldots,x_{k,n})$. Since by (A2),
$$
\del F(x_{1,n},\ldots,x_{k,n})\subset \del f_1(x_{1,n})\times\cdots\times
\del f_k(x_{k,n}),
$$
we have $\hat{x}^*_n=(x^*_{1,n},\ldots,x^*_{k,n})$ with
$x^*_{i,n}\in\del f_i(x_{i,n})$ for every $i=1,\ldots, k$.
We show that the sequences $((x_{i,n},x_{i,n}^*))_n\subset\del f_i$
satisfy all the requirements.
\smallbreak
From (a), we see that $x_{i,n}\to \xb$ for each $i$. Fix $\eps>0$ and
$j\in\{1,\ldots,k\}$.
By lower semicontinuity of the functions $f_i$, for $n$
sufficiently large it holds
\begin{equation}\label{sci1}
\sum_{i\ne j}f_i(\xb)\le \sum_{i\ne j}f_i(x_{i,n})+\eps,
\end{equation}
and by (b) above,
\begin{equation}\label{sci2}
F(\hat{x}_n)=\sum_{i}f_i(x_{i,n})\le F_\Delta(\xb,\ldots,\xb)=\sum_{i}f_i(\xb)+\eps.
\end{equation}
Combining \eqref{sci1} and \eqref{sci2}, we get that, for $n$
sufficiently large,
$$
\sum_{i}f_i(x_{i,n})\le \sum_{i\ne j}f_i(x_{i,n})+f_j(\xb)+2\eps
$$
hence, $f_j(x_{j,n})\le f_j(\xb)+2\eps$ eventually.
This shows that, for every $j$, $f_j(x_{j,n})\to f_j(\xb)$.
Finally, we have proved that $x_{i,n}\to_{f_i} \xb$ for every $i$.
\medbreak
We now show that $\limsup_n\|\sum x^*_{i,n}\|\le
\left |\nabla \left(\sum f_i\right)\right |(\xb)$.
Since $\hat{y}^*_n\in \del (nd_\Delta)(\hat{y}_n)$,
we have $\hat{y}^*_n\in \Delta^\perp$ and $\|\hat{y}^*_n\|\le n$.
On the other hand, $\|{\hat{x}_n^*}|\Delta\|=\|\sum x^*_{i,n}\|$ and
(c) implies
\begin{align*}
\limsup_n \|{\hat{x}_n^*}|\Delta\|\le
\limsup_n\, (\|\hat{x}^*_n+\hat{y}^*_n\|+\|{\hat{y}^*_n}|\Delta\|)
= \limsup_n\, (\|\hat{x}^*_n+\hat{y}^*_n\|
\le\left |\nabla \left(\sum f_i\right)\right |(\xb).
\end{align*}
So,
$\limsup_n \|\sum x^*_{i,n}\|\le\left |\nabla \left(\sum f_i\right)\right |(\xb)$.
\medbreak
It remains to show that $\diam (x_{1,n},\ldots,x_{k,n})\,\|x^*_{i,n}\|\to 0$.
Let $\eps>0$ be fixed. Since for each $i=1,\ldots, n$, $f_i(x_{i,n})\to f_i(\xb)$, we have
$|F(\hat{x}_n)-F_\Delta(\xb,\ldots,\xb)|\le \eps$ eventually,
so by (b) $d_\Delta(\hat{y}_n) \le \eps/n$ eventually, and also
$d_\Delta(\hat{x}_n) \le \eps/n$ eventually since $(\hat{y}_n)_n$
and $(\hat{x}_n)_n$ converge to the same point.
Let $\gamma=|\nabla \left(\sum f_i\right)|(\xb)$.
It follows from (c) that $\|\hat{x}^*_n+\hat{y}^*_n\|\le \gamma+\eps$
eventually, so $\|\hat{x}^*_n\|\le \gamma+\eps+n$ eventually.
Combining all of this, for $n$ sufficiently large it holds
$$d_\Delta(\hat{x}_n) \|\hat{x}^*_n\|\le \eps(\gamma+\eps)/n+\eps.$$
The conclusion follows directly from this last inequality.
This completes the proof of (R2).
\medbreak
Proceeding as in Theorem \ref{jetSlicestable}, it is easy to show that
(R3) follows from (R2); we omit the details.
\finpf

Rules (R2) and (R3) were considered and proved to be equivalent to (R1)
in \cite[Th\'eor\`eme 5.1]{Jul03}. The above approach through Wijsman stability
of subdifferentials is new.
More equivalent properties can be found in
\cite{Zhu98,Iof98,Las01}. See \cite{Iof12} for historical comments.

{\small

}
\end{document}